\documentclass[11pt]{article}
\usepackage{amsmath,amstext,amscd}
\usepackage{amssymb}
\usepackage{epsfig}
\usepackage{amsfonts,latexsym}

\newtheorem{thm}{Theorem}[section]
\newtheorem{lemma}[thm]{Lemma}
\newtheorem{cor}[thm]{Corollary}
\newtheorem{prop}[thm]{Proposition}

\newtheorem{Open questions}[thm]{Open questions}

\newtheorem{Problem}[thm]{Problem}
\newtheorem{Open Problem}[thm]{Open Problem}

\def\cal{\mathcal}

\def\Bbb{\mathbb}
\def\bar{\overline}

\newcommand{\F}{\mathbb{F}}
\def\Z{\Bbb{Z}}
\def\N{\Bbb{N}}

\def\ni{\noindent}

\def\Diam{\hbox{\rm Diam}}

\def\Cay{\hbox{\rm Cay}}

\def\SL{\hbox{\rm SL}}
\def\gcd{\hbox{\rm gcd}}
\def\PSL{\hbox{\rm PSL}}

\def\ms{\medskip}

\def\onto{{\kern3pt\to\kern-8pt\to\kern3pt}}

\def\<{\langle}
\def\>{\rangle}
\def\|{{\ |\ }}

 \def\AA{\cal A}
\def\BB{\cal B}
 \def\CC{\cal C}

\def\MM{\cal M}
 \def\PP{\cal P}
\def\NN{\cal N}

\newcommand{\set}[1]{\left\{#1\right\}}

\newcommand{\abs}[1]{\left|#1\right|}
\renewcommand{\ni}{\noindent}

\renewcommand{\ms}{\medskip}
\newcommand{\bs}{\bigskip}
\def\*{^{\star}}
\newcommand{\Proof}{\ni \emph{Proof. }}

\newcommand{\norm}[1]{|\!|#1|\!|}

\def\qed{{\ifhmode\unskip\nobreak\hfil\penalty50 \hskip1em \else\nobreak\fi
   \mbox{}\nobreak\hfil \rule{1ex}{1ex}%
   \parfillskip=0pt \finalhyphendemerits=0 \par}}

\begin{document}

\title{Navigating in the Cayley graphs \\ of $\SL_N(\Z)$ and $\SL_N(\F_p)$}

\author{T. R. Riley\thanks{Support from NSF grant  0404767 is gratefully acknowledged.}}

\date{October 2004}
\maketitle

\begin{abstract}
We give a non-deterministic algorithm that expresses elements of $\SL_N(\Z)$, for $N \geq 3$, as words in a finite set of generators, with the length of these words at most a constant times the word
metric.  We show that the non-deterministic time-complexity of the subtractive version of Euclid's algorithm for finding the greatest common divisor of $N \geq 3$ integers $a_1, \ldots, a_N$ is at most a constant times $N \log n$ where $n := \max\set {\abs{a_1}, \ldots, \abs{a_N}}$.  
This leads to an elementary proof that for $N \geq 3$ the word metric in $\SL_N(\Z)$
is biLipschitz equivalent to the logarithm of the matrix norm -- an instance of a theorem of
Mozes, Lubotzky and Raghunathan.  And we show constructively that there exists $K>0$ such that for all $N \geq 3$ and primes $p$, the diameter of the Cayley graph of $\SL_N(\F_p)$ with respect to the generating set $\set{ e_{ij} \mid i \neq j }$ is at most  $K N^2 \log p$.  \\ 
\medskip

\footnotesize{\ni \textbf{2000 Mathematics Subject
Classification: 20F05} \\ \ni \emph{Key words and phrases:} special linear, normal form, diameter, Cayley graph, Euclid's algorithm}
\end{abstract}

\section{Introduction} \label{intro}

This paper concerns expressing elements of $\SL_N(\Z)$ and $\SL_N(\F_p)$, for $N \geq 3$, as words in the generating set $\set{e_{ij} \mid i \neq j}$ consisting of the $N^2-N$ elementary matrices $e_{ij}$ that have $1$'s along the diagonal, the off-diagonal $ij$-entry $1$, and all other entries 0.  

What gets our study off the ground is an explicit means of writing powers ${e_{ij}}^m$ in $\SL_N(\Z)$, for $N \geq 3$, as products of $O(\log(1+ \abs{m}))$ matrices in $\set{e_{ij} \mid i \neq j}^{\pm 1}$.  This is explained in Section~\ref{compression section} and involves expressing  $m$ as a sum of Fibonacci numbers.  It is used (in Section~\ref{Euclid}) in a study of the non-deterministic time complexity of the subtractive Euclid's algorithm for finding the greatest common divisor of integers $a_1, \ldots, a_N$.  This differs from the standard Euclid's algorithm in that in each step one integer is added to or subtracted from another, rather than a remainder on division taken.  Yao and Knuth \cite{YK} proved that the \emph{average} number of steps to compute $\gcd(m,n)$ by the (deterministic) subtractive version of Euclid's algorithm, where $m$ is uniformly distributed in the range $1 \leq m \leq n$, is $6 \pi^{-2}(\ln n)^2 +O(\log n (\log  \log n)^2)$.  We show that the \emph{worst--case} non--deterministic complexity of Euclid's algorithm for computing the g.c.d.\ of $N \geq 3$ integers $(a^1, \ldots, a^N)$ is $O(N\log n)$, where $n:= \max\set{ \abs{a^1}, \ldots, \abs{a^N} }$.

\bs \ni \textbf{Theorem~\ref{subtractive alg}}
\emph{Suppose $(a^1, \ldots, a^N)$ is an $N$-tuple of integers, not all zero, and $N \geq 3$. Define $n:= \max\set{ \abs{a^1}, \ldots, \abs{a^N} }$.  There is a constant $K>0$, independent of $n$ and $N$, such that there is a sequence of no more than $K(N-1)(1+ \log n  )$  additions and subtractions of one entry from another, after which all but one entry in the $N$-tuple are zero.} 

\bs \ni The innovation is to use the compression techniques of Section~\ref{compression section} to accelerate repeated additions or subtractions of one entry to or from another.  By contrast,  the non--deterministic complexity is $\sim \!\! n$ in the case $N =2$  -- we supply a group theoretic proof of this, presumably well-known, result.  A vivid example is that it requires $\abs{n}$ steps (additions and subtractions) to convert $(1,n)$ to $(1,0)$, but $(1,n,0)$ can be reduced to $(1,0,0)$ in $O(\log \abs{n})$ steps. 

Then, in Section~\ref{MLR section}, we run our accelerated version of Euclid's algorithm of Section~\ref{Euclid} on the columns of a matrix $\MM$ in $\SL_N(\Z)$, in the course of reducing $\MM$ to the identity by row operations.  This leads to a new proof of an instance of a celebrated theorem  of Mozes, Lubotzky and Raghunathan \cite{LMR}, \cite{LMR2}, and we contribute information about the constants:

\bs \ni \textbf{Theorem~\ref{MLR}} \emph{Fix $N \geq 3$.  Let $\ell(\MM)$ denote the word length of $\MM \in \SL_N(\Z)$,
with respect to a fixed finite generating set. There
exist $C_1,C_2>0$ such that for all $\MM \in \SL_N(\Z)$
$$C_1\, \log \norm{\MM} \ \leq \ \ell(\MM) \ \leq \ C_2\, \log \norm{\MM}.$$
Moreover, if the generating set is $\set{e_{ij} \mid i \neq j}$ then $C_1$ is independent of $N$ and $C_2\leq C_3N^N$ for a constant $C_3>0$ that is independent of $N$.}

\bs
Our proof of Theorem~\ref{MLR}  is constructive (as are the proofs of the results in Sections~\ref{compression section} and \ref{Euclid} it appeals to) and amounts to an effective algorithm for finding a \emph{normal form} for $\SL_N(\Z)$ for $N \geq 3$ -- that is, for every $\MM \in \SL_N(\Z)$, a word $w_{\MM}$ on a fixed finite generating set and representing $\MM$. Equivalently, a normal form is a choice for all $\MM$ of path in the Cayley graph from the identity to $\MM$.  By homogeneity, it amounts to a means of navigating between any two vertices in the graph.  

Our normal form for $\SL_N(\Z)$ is \emph{of linearly bounded length}; that is, there exists $K >0$ such that for all $\MM \in \SL_N(\Z)$, the length of $w_{\MM}$ is at most $K$ times the length $\ell(\MM)$ of the shortest word that represents $\MM$.  This is because, the length of $w_{\MM}$ is at most $C_2 \log \norm{\MM}$, on account of its role in the proof of Theorem~\ref{MLR}, and $C_2 \log  \norm{M} \leq (C_2 /C_1) \ell(\MM)$.  

The author's original motivation for embarking on the work in this article was a potential application to the construction of van~Kampen diagrams to establish certain isoperimetric functions (concerning filling loops with discs): a long-standing claim of Thurston, originally quoted in \cite{Gersten} and repeated in \cite[$\S5.A_8$]{Gromov}, is that $\textup{SL}_N(\mathbb{Z})$ admits a quadratic isoperimetric function for all $N \geq 4$.  By contrast, Epstein and Thurston showed that the minimal isoperimetric function for $\SL_3(\Z)$ grows at least exponentially ~\cite[Chapter~10]{Epstein}.  By a theorem of Gromov~\cite[\S$5A_7$]{Gromov}, $\SL_3(\Z)$ admits an exponential isoperimetric function.   The author hopes the normal form will be of use towards proving Thurston's assertion and giving an elementary proof of Gromov's result. However there may be formidable obstacles; the geometry of the normal form has to be complicated in the following sense. For $N \geq 3$, no normal form for $\SL_N(\Z)$ of linearly bounded length can (either \emph{synchronously} or \emph{asynchronously}) \emph{fellow-travel}.  This result was proved by Epstein and Thurston~\cite[Chapter~10]{Epstein} to show that $\SL_N(\Z)$ is not automatic for $N \geq 3$; they use isoperimetric inequalities concerning filling $(N-2)$-spheres with $(N-1)$-balls (we mentioned the case $N=3$ above). 

\bs
Finally, in Section~\ref{diam section}, we apply similar techniques to $\SL_N(\F_p)$.  We find a normal form and prove the following result about Cayley graph diameter.

\bs \ni \textbf{Theorem~\ref{diam thm}}
\emph{There exists $C>0$ such that for all $N \geq 3$ and primes $p$, $$\Diam \ \Cay(\SL_N(\F_p), \set{ e_{ij} \mid i \neq j }) \ \leq \ C N^2 \log  p.$$} 

\bs
\ni We remark, for comparison, that a lower bound on the diameter of a constant times $( N^2 / \log N) \log  p$ follows from $\abs{ \SL_N(\F_p) } \sim p^{N^2 -1}$ because $\abs{\set{ e_{ij} \mid i \neq j }}=N^2-N$.

Define 
$$\mathcal{A}_N := \left(\begin{array}{ccccc}1 & 1 &  &  &  \\ & 1 &  &  &  \\ &  & 1 &  &  \\ &  &  & \ddots &  \\ &  &  &  & 1\end{array}\right), \ \ \  \mathcal{B}_N :=  \left(\begin{array}{llllll} \ \ \ 0 & 1 &  &  &  \\ & 0 & 1 &  &  \\ &  & 0 & \ddots &  \\ &  &  & \ddots & 1 \\\parbox{8mm}{$(-1)^{N-1}$} &  &  &  & 0\end{array}\right).$$
In Lemma~\ref{change} (which is due to M.~Kassabov) we show using elementary, constructive means that every $e_{ij}$ equals a word in ${\mathcal{A}_N}^{\pm 1}$ and ${\mathcal{B}_N}^{\pm 1}$ of length at most $10N$.  Applying this to Theorem~\ref{diam thm} we get:

\begin{cor} \label{diams cor}
There exists $C >0$ such that for all $N \geq 3$ and primes $p$,  $$\Diam \ \Cay(\SL_N(\F_p), \set{\mathcal{A}_N, \mathcal{B}_N}) \ \leq  \ C N^3 \log  p.$$  
\end{cor}

Lubotzky~\cite{Lubotzky} explains a non-constructive proof that given $N \geq 3$ and a generating set $X$ for $\SL_N(\Z)$, there exists   $K>0$,  that will depend on $N$ and potentially (see Problem~\ref{uniform diameters} below) on $X$, such that $\Diam \ \Cay(\SL_N(\F_p),X)  \ \leq \ K \log p$ for all primes $p$: since  $\SL_N(\Z)$ enjoys Property~(\emph{T}) for $N \geq 3$, the graphs $\set{\Cay(\SL_N(\Z/n\Z),X) \mid n \in \N }$ are a family of expanders, and the result follows.  This article supplies an elementary, constructive proof that avoids the big guns of Property (\emph{T}) and Selberg's Theorem.  The prior absence of such a proof is lamented on  of  \cite[page 102]{Lubotzky}.  

The argument above can be made quantitative as follows to yield a result that is weaker than that of Theorem~\ref{diam thm} in that it gives $N^3$ in place of the $N^2$ term in the estimate.  Kassabov~\cite{Kassabov}, extending methods of Shalom~\cite{Shalom2}, \cite{Shalom}, shows that
the Kazhdan constant for $\SL_N(\F_p)$ with respect to $\set{ e_{ij} \mid i\neq j}$ is at least $k:=(31\sqrt{N} +700)^{-1}$.    Define $\Gamma_{N,p}:= \Cay(\SL_N(\F_p), \set{ e_{ij} \mid i \neq j })$. The first non-zero eigenvalue $\lambda_1$ of the discrete Laplacian on $l_2(\Gamma_{N,p})$ is $1-s$ where $s$ is the spectral gap; $s \geq k^2/2$ by \cite{PZ}; and $$\Diam \, \Gamma_{N,p} \ \leq \ - (\log 2 \abs{ \SL_N(\F_p)})/(\log \lambda_1)$$ by \cite[Proposition~5.24]{LZ}.  In fact, any upper bound on diameter obtained this way is also an upper bound on \emph{mixing time of the random walk} on the Cayley graph, and so Theorem~\ref{diam thm} suggests that, as is often the case, mixing time and diameter differ for $\Gamma_{N,p}$.  

It is an open question\footnote{Added 31st January, 2005: this question has been answered in the affirmative by M.\ Kassabov and the author.}  \cite[Problem 8.1.3]{Lubotzky} whether the $N^3$ of Corollary~\ref{diams cor} can be improved to $N^2$.  Such a result would be best possible because $\abs{ \SL_N(\F_p) } \sim p^{N^2 -1}$.   Lubotzky, himself, gets close by proving with Babai and Kantor:

\begin{prop} \label{BKL prop} \emph{\cite{BKL}, \cite[Proposition 8.1.7]{Lubotzky}}
There exists $C >0$ such that for all $N \geq 3$ and primes $p$, there is a set $S$ of three generators for $\SL_N(\F_p)$ such that  $$\Diam \ \Cay(\SL_N(\F_p), S ) \ \leq  \ C N^2 \log  p.$$  In fact, $S$ can be taken to be $\set{ \AA_N, \BB_N, \CC_N}$ where $\AA_N$ and $\BB_N$ are defined above and $\CC_N := e_{12}{e_{21}}^{-1}e_{12}$.
\end{prop}

\ni The proof in \cite{Lubotzky} appeals to Selberg's Theorem, but the proof in \cite{BKL} is constructive and elementary save that ``unnatural'' generators of $\SL_2(\F_p) \leq \SL_N(\F_p)$ are used in place of $\AA_N$ and $\CC_N$.  An alternative route to Corollary~\ref{diams cor} is to apply Lemma~\ref{change} to Proposition~\ref{BKL prop}.

\bs

The following problems  provide a wider context for the study of diameters of Cayley graphs of $\SL_N(\F_p)$. 

\begin{Problem} \label{uniform T}  Fix $N \geq 3$. Does $\SL_N(\Z)$ enjoy uniform Property \emph{(}T\emph{)}?  
\end{Problem}

\begin{Problem} \label{big expander family} \emph{(}An Independence Problem for $\SL_N(\Z)$.\emph{)}
Fix $N \geq 3$.  Is $$\set{ \ \Cay(\SL_N(\Z) /H, X) \ \mid \ [\SL_N(\Z) : H ] < \infty, \ \langle X \rangle = \SL_N(\Z) \ }$$ a family of expanders?     
\end{Problem}

\begin{Problem} \label{uniform diameters} Fix $N \geq 2$.  Does there exist $K>0$ such that for all generating sets $X$ for $\SL_N(\Z)$ and all primes $p$  
$$\Diam \ \Cay(\SL_N(\F_p), X)  \ \leq \ K \log p?$$
\end{Problem}

\ni For fixed $N \geq 3$, an affirmative answer to Problem~\ref{uniform T} would imply an affirmative answer to \ref{big expander family}, and that, in turn, would imply an affirmative answer to \ref{uniform diameters}.  In the case $N =2$, the same implications apply between \ref{uniform diameters} and the following analogues of \ref{uniform T} and \ref{big expander family}:  does $\SL_2(\Z)$ enjoy uniform Property ($\tau$) with respect to congruence subgroups (``The Selberg Property'' \cite{LZ}), and is  $\set{  \Cay(\SL_2(\Z /m\Z), X)  \mid  \langle X \rangle = \SL_2(\Z)  }$ a family of expanders?    More details can be found in \cite{Lubotzky} and \cite{LZ}; groups in which the analogue of Problem~\ref{uniform T} has a negative answer are constructed in \cite{GZ}; the original (more general) \emph{independence problems} are in \cite{LW}; and a rare example of an independence result is due to Gamburd \cite{Gamburd} who (roughly speaking) finds a large class of generating sets $X$ for $\SL_2(\Z)$ and primes $p$ for which $\Cay(\SL_2(\F_p),X)$ forms a family of expanders.  

\ms 

We briefly mention related results for $\SL_N(\Z)$ and $\SL_N(\F_p)$ when $N =2$.  Property ($\tau$) is enjoyed by $\SL_2(\Z)$ as a consequence of Selberg's Theorem (see \cite{Lubotzky}, \cite{LZ}, \cite{Selberg}), and so for any \emph{fixed} finite generating set $X$ for $\SL_2(\Z)$,  we find $$\set{\Cay(\SL_N(\F_p),X) \mid p \textup{ prime} }$$ is a family of expanders. 
So there exists $K>0$ such that $$\Diam \ \Cay(\SL_2(\F_p),X)  \ \leq \ K \log p$$ for all primes $p$.  This proof (explained in \cite{Lubotzky}) is not constructive and neither is the only other known proof, which uses the \emph{circle method} for lifting elements of $\SL_2(\F_p)$ to elements of $\SL_2(\Z)$ with short word representations \cite{LPS}.  But Larsen~\cite{Larsen} has given an algorithm that produces word representations of length $O(\log p \log \log p)$.  In common with this article, representing powers such as ${e_{12}}^m$ by short words is key, and the subtractive version of Euclid's algorithm plays a role.    

Another constructive result is due to Gamburd and Shahshahani \cite{GS} and is in the direction of Problem~\ref{uniform diameters} in the case $N=2$. They give an algorithm that produces paths in Cayley graphs to prove the following uniform diameter bound: for all primes $p>2$, and for \emph{all} finite sets  $X$ of elements of $\PSL_2(\Z)$ such that $\langle X \rangle$ is a $p^2$-dense subgroup of $\PSL_2(\Z)$
$$\Diam \ \Cay(\PSL_2(\Z/p^n \Z), X )  \ = \ c  \log^d \abs{  \PSL_2(\Z/p^n\Z)},$$
where $d= \log_2 420$ and $c$ depends on $X$.  This has been recently improved by Dinai \cite{Dinai} who shows that for all $d>3$, there exists $c>0$ such that $\Diam \ \Cay(\SL_2(\Z/p^n \Z), X )  \leq c  \log^d \abs{  \SL_2(\Z/p^n\Z)}$ for \emph{all} generating sets $X$ for $\SL_2(\Z)$.

\bs
\ni \emph{Acknowledgements.}  I am grateful to Tsachik~Gelander, Martin~Kassabov and Alex~Lubotzky for explaining background to the subject of diameters of the Cayley graphs of $\SL_N(\F_p)$ to me, and to Karen Vogtmann for encouragement to investigate Thurston's claims about isoperimetric function for $\SL_N(\Z)$.  I additionally wish to thank Martin~Kassabov for providing Lemma~\ref{change}, improving a lemma in an earlier version of this article.

\section{Compressing powers ${e_{ij}}^m$} \label{compression section}

This section is devoted to proving the following result about representing powers ${e_{ij}}^m$ in $\SL_N(\Z)$ by words of length $O(\log(1 + \abs{m}))$.   

\begin{prop}\label{comp}
Suppose $N,m,i,j \in \Z$ with $N \geq 3$, with $1\leq i,j \leq N$, and with $i \neq j$. There exists a word $w_m \in \set{{e_{pq}}^{\pm1} \mid 1 \leq p, q \leq N, p \neq q}^{\star}$ such that $w_m = {e_{ij}}^m$ in $\SL_N(\Z)$ and $$\ell(w) \ \leq \  4+ 6\log_{\tau}(1+\abs{m}\sqrt{5}).$$
\end{prop}

It suffices to prove the result for $N=3, m >0, i=1$ and $j=3$, which we do by giving $w_m$ explicitly in the second of the two lemmas below.  The first lemma addresses the case where $m$ is a Fibonacci number (defined recursively by $F_0=0, \, F_1=1, \, F_{i+2}=F_{i+1}+F_i$), and will be superseded by the second lemma.  The detailed calculation in the proof of the first lemma is key to understanding the proof of the second.

\begin{lemma} \label{Fib powers} For  non-negative integers $n$, the words
\begin{tabbing}
 \ \ \ \ \= ${e_{23}}^{-1}(e_{23}e_{32})^{-n}{e_{13}}^{-1}(e_{23}e_{32})^{n}{e_{23}}^{-1}(e_{23}e_{32})^{-n} e_{13}(e_{23}e_{32})^n{e_{23}}^2,$ \= \  \  and \\  
\rule{0mm}{6mm}  \> ${e_{23}}^{-1}(e_{23}e_{32})^{-n}{e_{12}}^{-1}(e_{23}e_{32})^{n}{e_{23}}^{-1}(e_{23}e_{32})^{-n} e_{12}(e_{23}e_{32})^n{e_{23}}^2$ \> 
\end{tabbing}
equal ${e_{13}}^{F_{2n}}$ and ${e_{13}}^{F_{2n+1}}$, respectively, in $\SL_3(\Z)$.  
\end{lemma}

\Proof 
We multiply out the first of these words from right to left as follows.  The calculation for the second is very similar.   The notation for each step shown is $\mathcal{A} \xrightarrow{\mathcal{B}} \mathcal{B}\mathcal{A}$.   
\begin{tabbing}
AAA \= AAAAA \=  AAAAAAAAAAAAAA \= AAAAAA \=  \kill \\

$\left(\begin{array}{lll}1 & 0 & 0 \\0 & 1 & 0 \\0 & 0 & 1\end{array}\right)$   \ \  $\xrightarrow{ \ \ \parbox{8mm}{\scriptsize{${e_{23}}^2$}} }$
 \ \   $\left(\begin{array}{lll}1 & 0 & 0 \\0 & 1 & 2 \\0 & 0 & 1\end{array}\right)$   \> \> \>   $\xrightarrow{ \parbox{12mm}{\scriptsize{$(e_{23}e_{32})^n$}} }$  \>  $\left(\begin{array}{lll}1 & 0 & 0 \\0 & F_{2n+1} & F_{2n+3} \\0 & F_{2n} & F_{2n+2}\end{array}\right)$  \\

\rule{0mm}{12mm} \>  $\xrightarrow{ \ \parbox{8mm}{  \scriptsize{$ \ \ e_{13}$}} }$ \>  $\left(\begin{array}{lll}1 & F_{2n} & F_{2n+2} \\0 & F_{2n+1} & F_{2n+3} \\0 & F_{2n} & F_{2n+2}\end{array}\right)$ \> 
$\xrightarrow{ \parbox{12mm}{\scriptsize{$(e_{23}e_{32})^{-n}$}} }$ \>  $\left(\begin{array}{lll}1 & F_{2n} & F_{2n+2} \\0 & 1 & 2 \\0 & 0 & 1\end{array}\right)$ \\ 

\rule{0mm}{12mm} \> $\xrightarrow{ \ \parbox{8mm}{ \scriptsize{${e_{23}}^{-1}$}} }$ \>  $\left(\begin{array}{lll}1 & F_{2n} & F_{2n+2} \\0 & 1 & 1 \\0 & 0 & 1\end{array}\right)$ \>  $\xrightarrow{ \parbox{12mm}{\scriptsize{$(e_{23}e_{32})^n$}} }$ \>  $\left(\begin{array}{lll}1 & F_{2n} & F_{2n+2} \\0 & F_{2n+1} & F_{2n+2}  \\0 & F_{2n} & F_{2n+1}  \end{array}\right)$  \\ 

\rule{0mm}{12mm} \>  $\xrightarrow{ \ \parbox{8mm}{ \  \scriptsize{${e_{13}}^{-1}$}} }$  \>
$\left(\begin{array}{lll}1 & 0 & F_{2n} \\0 & F_{2n+1} & F_{2n+2}  \\0 & F_{2n} & F_{2n+1}  \end{array}\right)$ \>  $\xrightarrow{ \parbox{12mm}{\scriptsize{$(e_{23}e_{32})^{-n}$}} }$ \>   $\left(\begin{array}{lll}1 & 0 & F_{2n} \\0 & 1 & 1 \\0 & 0 & 1  \end{array}\right)$ \\

\rule{0mm}{12mm} \>  $\xrightarrow{ \ \parbox{8mm}{ \  \scriptsize{${e_{23}}^{-1}$}} }$  \>
$\left(\begin{array}{lll}1 & 0 & F_{2n} \\0 & 1 & 0  \\0 & 0 & 1 \end{array}\right)$ \> \>  
\end{tabbing}
\qed

\ms The following result can be proved by an easy induction. 

\bs \ni \textbf{Zeckendorf's Theorem} \cite{GKP},
\cite{Zeckendorf}. \emph{Every positive integer $m$ can be
expressed in a unique way as
\begin{equation} \label{Zeckendorf}
m= F_{k_1} + F_{k_2}+ \cdots + F_{k_r},
\end{equation}
with $k_1\geq 2$ and
$k_{j+1} - k_j \geq 2$ for all $1 \leq j < r$.}

\bs  
In fact, $F_{k_r}$ is the largest Fibonacci number no bigger than $m$, and
$F_{k_{r-1}}$ is the largest no bigger than $m-
F_{k_r}$, and so on.  Recall that $F_n = (\tau^n - (-\tau)^{-n})/\sqrt{5}$ for all $n$, where $\tau := (1+\sqrt 5)/2$, and so  
\begin{equation}
F_n \geq \frac{\tau^n -1}{\sqrt{5}}. \label{Fib ineq} 
\end{equation}  Thus, as $F_{k_r} \leq m$,  
\begin{eqnarray}
k_r & \leq & \log_{\tau} (1+ m \sqrt 5 ). \label{ineq 1}
\end{eqnarray}

\begin{lemma} \label{compress}
Suppose $m$ is a positive integer expressed as in \emph{(\ref{Zeckendorf})}.  Write $$m= (F_{\hat k_1} + F_{{\hat k}_2}+ \cdots + F_{{\hat k}_{\hat r}}) + (F_{{\bar k}_1} + F_{{\bar k}_2}+ \cdots + F_{{\bar{k}}_{\bar{r}}})$$ where ${\hat k}_1<  \ldots < {\hat k}_{\hat r}$ are the even numbers amongst $k_1, \ldots, k_r$ and ${\bar k}_1<  \ldots < {\bar k}_{\bar r}$ are the odd numbers. 
Let $n$ be the integer such that either $2n=k_r$ or $2n+1=k_r$.  Let $u_m$ be the word
$$a_{n}b_{n} (e_{23} e_{32})\ldots a_2b_2(e_{23}e_{32}) a_{1}b_{1}(e_{23}e_{32})$$
in which $a_i = e_{13}$ if $2i \in \lbrace {\hat k}_1, \ldots, {\hat k}_{\hat r} \rbrace$ and is the empty string otherwise, and $b_i=e_{12}$ if $2i+1 \in \set{ {\bar k}_1, \ldots, {\bar k}_{\bar r} }$ and is the empty string otherwise. 
Let  $v_m$ be the word obtained from $u_m$ by replacing every $e_{12}$ and $e_{13}$ by ${e_{12}}^{-1}$ and ${e_{13}}^{-1}$, respectively.
 Define $$w_m := {e_{23}}^{-1}(e_{23}e_{32})^{-n}v_m{e_{23}}^{-1}(e_{23}e_{32})^{-n} u_m {e_{23}}^2.$$
Then $w_m$ equals ${e_{13}}^{m}$ in $\SL_3(\Z)$ 
and has length
\begin{eqnarray}
\ell(w_m) \  \leq \ 4+ 6 \log_{\tau}(1+m \sqrt{5}). \label{log ineq}
\end{eqnarray} 
\end{lemma}

\Proof
Lemma~\ref{Fib powers} is a special case of this lemma: when $m=F_{2n}$ we find $u_m = e_{13} (e_{23}e_{32})^n$ and $v_m = {e_{13}}^{-1} (e_{23}e_{32})^n$, and when $m=F_{2n+1}$ we find $u_m = e_{12} (e_{23}e_{32})^n$ and $v_m = {e_{12}}^{-1} (e_{23}e_{32})^n$.  Multiply out $w_m$ from right to left, as follows, using a more general and concise version of the calculation used to establish Lemma~\ref{Fib powers}. All the sums are over $i= 1,\ldots, r$.  

\begin{tabbing}
AAA \= AAA \= AAAAAAAAA \=  AAAAAAAAAAAAAA  \kill \\

\> $\left(\begin{array}{ccc}1 & 0 & 0 \\0 & 1 & 2 \\0 & 0 & 1\end{array}\right)$   \ \ \   $\xrightarrow{ \parbox{12mm}{\scriptsize{ \ \ \ \ \ $u_n$}} }$  \ \ \  $\left(\begin{array}{ccc}1 & \sum{F_{k_i}} & \sum{F_{k_i+2}} \\0 & F_{2n+1} & F_{2n+3} \\0 & F_{2n} & F_{2n+2}\end{array}\right)$  \> \> \\  

\> \rule{0mm}{12mm} \>  $\xrightarrow{ \parbox{20mm}{\scriptsize{ \ \ \ \ $(e_{23}e_{32})^{-n}$}} }$ \>  $\left(\begin{array}{ccc}1 & \sum{F_{k_i}} & \sum{F_{k_i+2}} \\0 & 1 & 2 \\0 & 0 & 1\end{array}\right)$  \\ 

\> \rule{0mm}{12mm} \>  $\xrightarrow{ \parbox{20mm}{ \ \ \ \ \scriptsize{$v_m{e_{23}}^{-1}$}} }$ \>  $\left(\begin{array}{ccc}1 & 0 & \sum({F_{k_i+2}}-{F_{k_i+1}}) \\0 & F_{2n+1} & F_{2n+2} \\0 & F_{2n} & F_{2n+1}\end{array}\right)$ \\

\> \rule{0mm}{12mm} \>  $\xrightarrow{ \parbox{20mm}{\scriptsize{${e_{23}}^{-1}(e_{23}e_{32})^{-n}$}} }$ \>  $\left(\begin{array}{ccc}1 & 0 & \sum{F_{k_i}} \\0 & 1 & 0 \\0 & 0 & 1 \end{array}\right).$ 
\end{tabbing}
The length of $w_m$ is $4+ 8n+2r$, from which we get (\ref{log ineq}) by using (\ref{ineq 1}), $r \leq k_r$ and $n \leq k_r /2$. 
\qed

\section{Accelerating the subtractive version of Euclid's algorithm}\label{Euclid}

The \emph{subtractive} version of Euclid's algorithm for finding the greatest common divisor of an $N$-tuple of integers differs from the standard Euclid's algorithm in that at  each \emph{step} an addition or subtraction is made rather than a remainder taken.  
That is, in one step an $N$-tuple $(a_{i+1}^1, \ldots, a_{i+1}^N)$ is produced from the previous $k$-tuple $(a_i^1, \ldots, a_i^N)$, as follows.     Take $p$ and $q$  so that $a_i^p$ and $a_i^q$ have the greatest and second greatest absolute values  amongst  $a_i^1, \ldots, a_i^N$.  (To resolve dead-heats, take $p$ minimal, and then take $q$ minimal amongst the remaining indices.)  Define $a_{i+1}^p := a_{i}^p \pm a_{i}^p$, with the sign chosen so that  $\abs{a_{i+1}^p} < \abs{a_{i}^p}$, and define $a_{i+1}^j := a_{i}^j$ for all $j \neq p$.  Stop when all but one entry is zero and output the absolute value of that entry.  

For example, in 6 steps the algorithm gives $\gcd(-32,8,-12) = 4$: 
 $$(-32,8,-12)  \ \mapsto \ (-20,8,-12) \ \mapsto \ (-8,8,-12) \ \mapsto \ (-8,8,-4) $$ 
 $$ \mapsto \  (0,8,-4) \  \mapsto \ (0,4,-4) \ \mapsto \ (0,0,-4). $$

There is a non--deterministic version of this algorithm in which obtaining $(a_{i+1}^1, \ldots, a_{i+1}^N)$ from $(a_i^1, \ldots, a_i^N)$ by adding one entry to another or by subtracting one entry from another constitutes a \emph{step}.  Again, the algorithm stops when all but one entry is zero, and the output is the absolute value of that entry.   

Yao and Knuth~\cite{YK} proved that the average number of steps to compute $\gcd(m,n)$ by the (deterministic) subtractive version of Euclid's algorithm, where $m$ is uniformly distributed in the range $1 \leq m \leq n$, is $6 \pi^{-2}(\ln n)^2 +O(\log n (\log  \log n)^2)$.  We will show that the worst--case non--deterministic complexity of Euclid's algorithm for computing the $\gcd$ of $N \geq 3$ integers $(a^1, \ldots, a^N)$ is $O(\log n)$, where $n:= \max\set{ \abs{a^1}, \ldots, \abs{a^N} }$.  (In particular, the greatest common divisor of two integers $(a^1,a^2)$ can be  calculated non-deterministically in $O(\log(n))$ steps by starting with $(a^1,a^2,0)$.) That is, we prove: 

\begin{thm} \label{subtractive alg}
Suppose $(a^1, \ldots, a^N)$ is an $N$-tuple of integers, not all zero, and $N \geq 3$. Define $n:= \max\set{ \abs{a^1}, \ldots, \abs{a^N} }$.  There is a constant $K>0$, independent of $n$ and $N$, such that there is a sequence of no more than $K(N-1)(1+ \log n  )$ additions and subtractions of one entry to or from another, after which all but one entry in the $N$-tuple are zero.    
\end{thm}

\Proof
First consider running the standard Euclid's algorithm on the first two entries $a_0:=a^1$ and $b_0:=a^2$ in the $N$-tuple.  This proceeds via a sequence $(a_i,b_i)$ of pairs of integers finishing with a pair $(a_k,b_k)$ one of which is zero.  The pair $(a_{i+1},b_{i+1})$ is obtained from $(a_i,b_i)$ by replacing the entry with the larger absolute value by the remainder on division by the other.  So for all $i=0, \ldots, k-1$ there is some integer $q_i$ such that either ($a_{i+1} = a_i \pm q_i b_i$ and $b_{i+1} = b_i$), or ($a_{i+1} = a_i$ and $b_{i+1} = b_i \pm q_i a_i$).

It takes the standard subtractive algorithm $q_i$ steps to get from $(a_i,b_i)$ to $(a_{i+1},b_{i+1})$.  But, as $N \geq 3$, Proposition~\ref{comp} gives us a word $w^{q_i}$ that has length at most $4 + 6 \log_{\tau} (1+q_i \sqrt 5))$ and that, reading right-to-left, describes a sequence of steps with the same effect.   (The step described by the letter ${e_{pq}}^{\pm 1}$ corresponds to left-multiplying the transpose of the $N$-tuple. The entries $a^{3},\ldots, a^N$ in the $N$-tuple may be disturbed in the course of these steps, but are recovered.)  

Define $c_i := \max \set{ \abs{a_i}, \abs{b_i}}$.  Then $q_i \leq {c_i} / c_{i+1}$ and $q_0 \leq n$.  So it is possible to get from $(a_0, b_0)$ to $(a_k, b_k)$ in $S$ steps where
$$ S \ \leq \ \sum_{i=0}^{k-1} \left(4+ 6\log_{\tau}\left(1+ \frac{c_i}{c_{i+1}} \sqrt{5}\right) \right).$$
But, as $c_i / c_{i+1} \geq 1$ for all $0 \leq i \leq k-1$, and $c_0 /c_k \leq n$, this is at most
\begin{eqnarray} 
S & \leq & 4k + 6 \sum_{i=0}^{k-1}  \log_{\tau}\left(\frac{c_i}{c_{i+1}} \left(1+ \sqrt{5}\right) \right) \nonumber \\ & \leq & 4k+6\left(\log_{\tau} n + k \log_{\tau} \left(1+ \sqrt 5\right)\right). \label{getting close}
\end{eqnarray}
Now $n \geq F_{k+1}$ by an easy induction.  So by inequality~(\ref{Fib ineq}) of Section~\ref{compression section} $$n  \ \geq \ \frac{\tau^{k+1} -1}{\sqrt 5},$$ and thus $k \leq -1+ \log_{\tau}(1+n \sqrt 5)$. This inequality  together with (\ref{getting close}) shows there exists $K >0$ such that $S \leq K +K \log n$.  

Obtain the bound claimed in the theorem by next arguing as above for $a^3$ and whichever or the first and second entries in the $N$-tuple is now non-zero, and then similarly for $a^4$, and so on, until finally for $a^N$.    
 \qed

\bs 
The proof above can be developed into a deterministic algorithm to calculate $\gcd(a^1, \ldots, a^N)$.  What are needed are the $q_i$ together with the words $w^{q_i}$ of Lemma~\ref{compress}.  But those $w^{q_i}$ are built using the expression for $q_i$ of Zeckendorf's Theorem.  Whilst is not hard to write routines to supply the $q_i$ and the expressions as per Zeckendorf's Theorem, it is not clear that producing a deterministic algorithm to calculate $\gcd(a^1, \ldots, a^N)$ in this manner has any computational advantages.

\bs

Theorem~\ref{subtractive alg} fails when $N=2$ (it is likely the following results are well known, but we include them for completeness and for the contrast):

\begin{prop}
To convert $(1,n)$ to $(\pm 1,0)$ or $(0, \pm 1)$ by successively subtracting one entry from, or adding one entry to, the other, requires $n$ steps.      
\end{prop}

\Proof
The number of steps required is at least the distance from ${e_{21}}^n$ to the identity in the word metric on $\SL_2(\Z)$ with respect to the generating set $e_{12}, e_{21}$.  This is because reading a word $w$ that represents ${e_{21}}^n$  from right to left would give a sequence of steps that transforms $(1,0)^t$ to $(1,n)^t$.   

But such a word $w$ descends to a word $\hat{w}$ in the images ${\mbox{$\hat{e}_{12}$}}^{\pm 1}$, ${\mbox{$\hat{e}_{21}$}}^{\pm 1}$ of  ${e_{12}}^{\pm 1}$, ${e_{21}}^{\pm 1}$ under the natural map $\SL_2(\Z) \onto \PSL_2(\Z)$.  And $\PSL_2(\Z) \cong (\Z / 2\Z) \ast (\Z / 3\Z)$, presented by $\langle \hat{s}, \hat{t} \mid \hat{s}^2 , \hat{t}^3\rangle$, where $$s \  = \ \left(\begin{array}{cc}0 & -1 \\1 & 0\end{array}\right) \  \ \textup{and} \  \  t \  = \ \left(\begin{array}{cc}1 & 1 \\ -1 & 0\end{array}\right).$$ Now, $st=e_{21}$ and so $(st)^n = {e_{21}}^n$.   And $(\hat{s}\hat{t})^n$ is of minimal length amongst all words in $\set{{\hat{s}}^{\pm 1}, {\hat{t}}^{\pm 1}}^{\star}$ that represent  $\mbox{$\hat{e}_{21}$}^n$ in the free product $(\Z / 2\Z) \ast (\Z / 3\Z)$.  So the minimal length of words in $\set{{e_{12}}^{\pm 1}, {e_{21}}^{\pm 1}}^{\star}$ that equal $\mbox{${e}_{21}$}^n$ in $\SL_2(\Z)$ is $n$.
 \qed
 
 \begin{cor}
 The \emph{(}worst case\emph{)} non-deterministic time complexity of the subtractive version of Euclid's algorithm for finding the greatest common divisor of two integers $a,b$ with  $n:= \max\set{\abs{a},\abs{b}}$ is between $n$ and $2n$. 
 \end{cor}

In the next section we will need the following more technical result that is proved in the same way as Theorem~\ref{subtractive alg}.

\ms \ni \textbf{Theorem~\ref{subtractive alg}$'$}
\emph{Suppose $(a^1, \ldots, a^N)$ is an $N$-tuple of integers, not all zero, where $N \geq 3$, and suppose $1 \leq k \leq N$. Define $n:= \max\set{ \abs{a^{N-k+1}}, \ldots, \abs{a^N} }$.  There is a constant $K>0$, independent of $k,n$ and $N$, such that there is a sequence of no more than $(k-1)K(1+ \log n)$ steps after which 
the first $N-k$ entries in the $N$-tuple are unchanged, all but one of the remaining entries in the $N$-tuple are zero, and that remaining entry is $\pm \gcd(a^{N-k+1}, \ldots, a^N)$.}

\section{The Mozes-Lubotzky-Raghunathan Theorem}  \label{MLR section}

In this section we give an elementary proof of the following result which is an instance of a theorem  of Mozes, Lubotzky and Raghunathan on irreducible lattices in semi-simple Lie groups of rank at least 2.  In \cite{LMR} they proved the case addressed below before generalising it to lattices in other  Lie groups in \cite{LMR2}.  We add information about the constants.  (For a matrix $\MM$ with real entries, $\norm{\MM}$ denotes the sup-norm, the maximum of the absolute values of the entries.)

\begin{thm} \label{MLR}  Fix $N \geq 3$.  Let $\ell(\MM)$ denote the word length of $\MM \in \SL_N(\Z)$,
with respect to a fixed finite generating set. There
exist $C_1,C_2>0$ such that for all $\MM \in \SL_N(\Z)$
$$C_1\, \log \norm{\MM} \ \leq \ \ell(\MM) \ \leq \ C_2\, \log \norm{\MM}.$$
Moreover, if the generating set is $\set{e_{ij} \mid i \neq j}$ then $C_1$ is independent of $N$ and $C_2\leq C_3N^N$ for a constant $C_3>0$ that is independent of $N$.   
\end{thm}

\Proof  One easily checks that if the first part of the theorem holds for one finite
generating set for $\SL_N(\Z)$ then it holds for all.  We will work with the generating set $\set{e_{ij} \mid i
 \neq j}$.

The first inequality is straightforward.  The sup-norm of a matrix that is the product of $n$ matrices in
$\set{e_{ij} \mid i \neq j}$ is at most $F_n$, and $F_n$ grows
exponentially with $n$.

The second inequality will take more work.  Suppose $\MM \in \SL_N(\Z)$.  Below, is a (well-known) procedure for reducing $\MM$ to the identity by row operations.  Each row operation corresponds to left-multiplication by some ${e_{ij}}^{\pm 1}$ and so a word $w \in \set{ {e_{ij}}^{\pm 1} \mid i \neq j }^{\star}$ that equals $\MM$ in $\SL_N(\Z)$ can be extracted.  
\begin{itemize}
\item[($1$)] Convert $\MM$ to an upper triangular matrix whose diagonal entries are all $\pm 1$, as follows. 
\begin{itemize}
\item[($1_1$)] Run Euclid's algorithm on the first column.  This will leave all entries zero except one that is $\pm 1$, because $\det \MM =1$.  Let $i_1$ be the row containing the non-zero entry in the first column.  If $i_1 \neq 1$ then premultiply by $e_{1i_1}{e_{i_11}}^{-1}e_{1i_1}$, which reverses the signs of the entries in row 1 and then interchanges rows $1$ and $i_1$.  
\item[($1_2$)] Run Euclid's algorithm on the entries in rows $2$ to $N$ of second column, leaving all  zero except one that is $\pm 1$ and lies in row $i_2$.  If $i_2 \neq 2$ then premultiply by $e_{2i_2}{e_{i_22}}^{-1}e_{2i_2}$.   
\item[\vdots]
\item[($1_{N\!-\!1}$)] Run Euclid's algorithm on the entries in rows $N-1$ and $N$ of the $(N-1)$-st column,  to make one entry $0$ and the other $\pm 1$. Then, if necessary, premultiply by $e_{N-1,N}{e_{N,N-1}}^{-1}e_{N-1,N}$ to get an upper triangular matrix.             
\end{itemize}
\item[($2$)] Get a matrix $(m_{ij})$ for which all the entries on the diagonal are 1 by premultipling by at most $N/2$ matrices $(e_{ij}{e_{ji}}^{-1}e_{ij})^2$ that reverse the signs of all the entries in rows $i$ and $j$.
\item[($3$)] Clear all the above--diagonal entries in $(m_{ij})$, one column at a time, as follows. 
\begin{itemize}
\item[($3_2$)] Premultiply by ${e_{12}}^{-m_{12}}$.   
\item[($3_3$)] Premultiply by ${e_{13}}^{-m_{13}}{e_{23}}^{-m_{23}}$.   
\item[\vdots] 
\item[($3_N$)]  Premultiply by ${e_{1,N}}^{-m_{1,N}} \ldots {e_{N-1,N}}^{-m_{N-1,N}}$.    
\end{itemize}
\end{itemize}

As it stands, the number of  ${e_{ij}}^{\pm 1}$ used in the procedure above may wildly exceed $\log \norm{\MM}$ on account of steps ($1$) and ($3$).  However, we can \emph{accelerate} ($1_1$)--($1_{N\!-\!1}$) as per Theorem~\ref{subtractive alg}$'$.  Define $\MM_0:=\MM$.  Performing ($1_1$) then takes at most $k_1:=(N-1)K(1+ \log \norm{\MM})$ steps and leaves a matrix $\MM_1$ such that $\norm{\MM_1} \leq \norm{\MM_0} F_{k_1+2}$.  
And, proceeding inductively, ($1_i$) costs at most $k_i:= (N-i)K(1+ \log \norm{\MM_{i-1}})$ steps and leaves a matrix $\MM_i$ with $\norm{\MM_i} \leq \norm{\MM_{i-1}} F_{k_i+2}$. 

So there is a constant $C>0$ such that for all $1\leq i \leq N-1$, 
\begin{eqnarray*}
k_i & \leq & KN(1+\log \norm{\MM_{i-1}}) \ \  \textup{and}  \\ 
\log{\norm{\MM_i}} & \leq & Ck_i+\log\norm{\MM_{i-1}}.  
\end{eqnarray*}
These inequalites and induction can be used to establish that for $1\leq i \leq N$
\begin{eqnarray}
\log \norm{\MM_{i-1}} & \leq &  (1 + \log \norm{\MM} ) ( 1+ CKN)^{i-1} -1, \label{ln}  
\end{eqnarray}
and for all $1\leq i \leq N-1$
\begin{eqnarray*}
k_i & \leq & KN (1 + \log \norm{\MM} )(1+CKN)^{i-1}. 
\end{eqnarray*}
So there is a constant $C' >0$, independent of $\MM$ and $N$, such that the contribution of ($1$) to $\ell(w)$ is at most $\sum_{i=1}^{N-1}k_i \leq C' N^{N-1} \log  \norm{\MM}$.

The contribution of step ($2$) to $\ell(w)$ is at most $3N$.  To assess the contribution of step ($3$), first note that $\abs{m_{ij}} \leq \norm{\MM_i}$ for all $i,j$ because in the course of step ($1$), row $i$ is not disturbed after ($1_i$). 
So, by inequality (\ref{ln}) and by compressing  each of the $N(N-1)/2$ terms ${e_{ij}}^{-m_{ij}}$ as per Proposition~\ref{comp}, we see that the effect of step ($3$) can be achieved whilst contributing at most 
\begin{equation} 
C'+ C' (1+\log \norm{\MM}) \sum_{i=1}^{N-1} (N-i)(1+CKN)^i \label{nasty sum}  
\end{equation}
to $\ell(w)$, for some constant $C'>0$ independent of $\MM$ and $N$.  But the summation term in (\ref{nasty sum}) is at most a constant times $1+N+N^2+ \cdots + N^N$, which is at most $4N^N$ as $N > 2$.  The outstanding claims of the theorem then follow.      
\qed

\section{The diameter of $\SL_N(\F_p)$.}  \label{diam section}

We adopt the notation $\alpha^{\beta}:= \beta^{-1} \alpha \beta$ and $[\alpha,\beta]= \alpha^{-1}\beta^{-1} \alpha \beta$.  

\begin{thm} \label{diam thm}
There exists $C>0$  such that for all $N \geq 3$ and primes $p$, $$\Diam \ \Cay(\SL_N(\F_p), \set{ e_{ij} \mid i \neq j }) \ \leq \ C N^2 \log  p.$$
\end{thm}

\Proof
Suppose $\MM \in \SL_N(\F_p)$.  We reduce $\MM$ to the identity matrix by successively premultiplying by matrices in $\set{ e_{ij} \mid i \neq j }$ in a similar manner to that used to prove Theorem~\ref{MLR}.      

First lift the entries in the first column of $\MM$ to $[0,p-1]$ and run the accelerated version of the subtractive version of Euclid's algorithm on the first column of the matrix.  Then swap two rows (changing the sign of one) to move the non-zero entry to the first row.  Next run the accelerated version of the subtractive version of Euclid's algorithm on the lift to $[0,p-1]$ of all but the first entry of second column, and move the non-zero entry to place 2,2 in the matrix.  Continue similarly through all the columns.  
By Theorem~\ref{subtractive alg}$'$, the cost is at most a constant times 
$$ \log p \sum_{j=1}^{N-1} (N-j) \ < \  N^2   \log p.$$  

We now have an upper triangular matrix $(m_{ij})$ such that every diagonal entry is non-zero.  As $p$ is prime and the diagonal entries in the matrix are non-zero, we can clear all the $N(N-1)/2$ entries above the diagonal by premultiplying by matrices of the form ${e_{ij}}^{-m_{ij}}$ where $j>i$.  By Proposition~\ref{comp} the effect of premultiplying by ${e_{ij}}^{-m_{ij}}$ can be achieved by premultiplying by a sequence of at most a constant times $\log p$ matrices in $\set{{e_{ij}}^{\pm 1} \mid i \neq j}$.  So we can reduce $\MM$ to a diagonal matrix $\mathcal{D} = \textup{diag}(a_1, \ldots, a_N)$ with total cost at most a constant times $N^2 \log p$.   

As $a_1$ and $a_1a_2$ are invertible in $\F_p$, we  can convert $\mathcal{D}$ to $\textup{diag}(1,a_1a_2, a_3, \ldots, a_N)$ by premultiplying by ${e_{12}}^{\pm 1}, {e_{21}}^{\pm 1}$ as follows, 
\begin{tabbing}
AAAAAAAAAA \= AAAAAA \=  AAAAAAAAAAA \= AAAAAAA  \=   \kill \\
$\left(\begin{array}{ccc}a_1 & 0 & \\0 & a_2  \\ & & \ddots  \end{array}\right)$   \>   $\xrightarrow{ \parbox{15mm}{\scriptsize{ \ \ \ \ \ $e_{12}{e_{21}}^{-1}e_{12}$}} }$  \>  $\left(\begin{array}{ccc}0 & a_2 & \\-a_1 & 0 & \\ & & \ddots \end{array}\right)$ \>  $\xrightarrow{ \parbox{17mm}{\scriptsize{ \ ${e_{12}}^{-{a_1}^{-1}}$}}}$ \>  $\left(\begin{array}{ccc}1 & a_2 &  \\ -a_1 & 0 & \\ & & \ddots \end{array}\right)$  \\

\rule{0mm}{12mm}  \>  $\xrightarrow{ \parbox{15mm}{\scriptsize{ \ ${e_{21}}^{a_1}$}}}$ \>  $\left(\begin{array}{ccc}1 & a_2 & \\ 0 & a_1a_2 & \\ & & 0 \end{array}\right)$  \>  $\xrightarrow{ \parbox{17mm}{\scriptsize{ \ ${e_{12}}^{- a_2(a_1a_2)^{-1}}$}}}$ \>  $\left(\begin{array}{ccc}1 & 0 & \\ 0 & a_1a_2 & \\ & & \ddots \end{array}\right).$
\end{tabbing}
Using Proposition~\ref{comp} the same effect can achieved by pre-multiplying by at most a constant times $\log p$ matrices in $\set{{e_{ij}}^{\pm 1} \mid i \neq j}$.  Applying this same  process to the second and third rows, and then  the third and fourth, and so on we reduce the matrix to the identity, at a total cost of at most a constant times $N \log p$.     
\qed

\bs
To deduce Corollary~\ref{diams cor} we use the following lemma, due to M.~Kassabov, concerning the matrices $\AA_N$ and $\BB_N$ given in Section~\ref{intro}.

\begin{lemma} \label{change} For all $1 \leq i,j \leq N$ with $i \neq j$ it is possible to express $e_{ij}$ as a word  in ${\mathcal{A}_N}^{\pm 1}$ and ${\mathcal{B}_N}^{\pm 1}$  of length at most $10N$. 
\end{lemma}

\Proof  We will drop the subscripts from $\mathcal{A}_N$ and $\mathcal{B}_N$. 
For all $i \neq j$ we have  ${e_{i,j}}^{\BB} = {e_{i+1,j+1}}^{\pm 1}$ where the indices are  in $\set{1, \ldots, N}$ and are taken modulo $N$.   So it suffices to express all $e_{1,1+k}$, for $1 \leq k \leq N-1$, as words in ${\mathcal{A}}^{\pm 1}$ and ${\mathcal{B}}^{\pm 1}$  of length at most $8 N$.    

For $ k =2, \ldots, n$ define $\PP_k:= e_{12}e_{23} \ldots e_{k-1,k}$.  Then $$\PP_k = \AA \AA^\BB\cdots \AA^{\BB^{k-2}} = \AA\BB^{-1}\AA\BB^{-1} \ldots \AA\BB^{-1} \AA{\BB^{k-2}}.$$ 
Now $\NN_k  \ := \  \PP_k  {\PP_{k-1}}^{-1}$
which, due to cancellations, equals a word of length $4k-7$ in $\AA^{\pm 1}$ and $\BB^{\pm 1}$, and is
\begin{eqnarray*} 
\left(\begin{array}{cccccc}
   \!\!1 & \!\!1 & \!\!\hdots & \!\!1 &\!\!1 & \\ 
   & \!\!1 & \!\!\ddots & \!\!\vdots & \!\!\vdots  &  \\ 
     &  & \!\!\ddots & \!\!1 & \!\!1 &  \\ 
   &  &  & \!\!1 & \!\!1 &  \\
   &  &  &  & \!\!1 &  \\
    &  &  &  & & \!\!\!\parbox{8mm}{$\textit{Id}_{n-k}$}
\end{array}\right)
\left(\begin{array}{rrrrrr}
1 & \!\!\!-1 &  &  &  &  \\ 
  & 1 & \!\!\!\ddots &  &  & \\ 
  &  & \!\!\!\ddots & \!\!\!-1 &   &  \\ 
   &  &  & 1 &  &  \\ 
 &  &  &  & \!\!\!1 &  \\ 
 &  &  &  &  & \!\!\!\parbox{8mm}{$\textit{Id}_{n-k}$}
\end{array}\right)    = 
\left(\begin{array}{cccccc}
1 &  &  &  & \!\!1 &  \\ 
  & \!\!1 &  &  & \!\!1 & \\ 
  &  & \!\!\ddots &  & \!\!\vdots  &  \\ 
   &  &  & \!\!1 & \!\!1 &  \\ 
 &  &  &  & \!\!1 &  \\ 
 &  &  &  &  & \!\!\parbox{8mm}{$\textit{Id}_{n-k}$}
\end{array}\right).
\end{eqnarray*}
For $k =3, \ldots, n$, we calculate that $\NN_k (\BB^{-1}\NN_{k-1}\BB)^{-1}$ is
\begin{eqnarray*} 
 \left(\begin{array}{cccccc}
\!\!1 &  &  &  & \!\!1 &  \\ 
  & \!\!1 &  &  & \!\!1 & \\ 
  &  & \!\!\ddots &  & \!\!\vdots  &  \\ 
   &  &  & \!\!1 & \!\!1 &  \\ 
 &  &  &  & \!\!1 &  \\ 
 &  &  &  &  & \!\!\!\parbox{8mm}{$\textit{Id}_{n-k}$}
\end{array}\right)
\left(\begin{array}{rrrrrr}
\!\!1 &  &  &  & 0 &  \\ 
  & \!\!1 &  &  & \!\!\!-1 & \\ 
  &  & \!\!\ddots &  & \vdots  &  \\ 
   &  &  & \!\!1 & \!\!\!-1 &  \\ 
 &  &  &  & \!\!1 &  \\ 
 &  &  &  &  & \!\!\!\parbox{8mm}{$\textit{Id}_{n-k}$}
\end{array}\right) \ = \ e_{1k}.
\end{eqnarray*}
So $e_{1k}$ can be expressed as a word of length $8k -16$ in $\AA^{\pm 1},\BB^{\pm 1}$.
\qed

\bibliographystyle{plain}
\bibliography{bibli}

\ni  \textsc{Tim R.\ Riley} \rule{0mm}{6mm} \\ Mathematics
Department, 10 Hillhouse Avenue, P.O. Box 208283, New Haven, CT 06520-8283, USA \\
\texttt{tim.riley@yale.edu, \
http:/\!/www.math.yale.edu/users/riley/ }

\end{document}